\documentclass{amsart}
\usepackage{xcolor}

\usepackage{graphicx}

\usepackage[utf8]{inputenc}
\usepackage{amsmath, amssymb}
\usepackage{amsfonts}
\usepackage{amsthm}
\usepackage{enumerate}
\usepackage{hyperref}
\usepackage{cleveref}
\usepackage[top=1in, bottom=1.3in, left=1.3in, right=1in]{geometry}
\newtheorem{theorem}{Theorem}[section]
\newtheorem{corollary}[theorem]{Corollary}
\newtheorem{lemma}[theorem]{Lemma}
\newtheorem{proposition}[theorem]{Proposition}

\theoremstyle{definition}
\newtheorem{remark}[theorem]{Remark}

\DeclareMathOperator{\im}{Im\,}

\newcommand{\Cm}{\mathbb{C}}
\newcommand{\R}{\mathbb{R}}
\newcommand{\Z}{\mathbb{Z}}

\newcommand{\N}{\mathbb{N}}

\newcommand{\T}{\mathbb{T}}
\newcommand{\C}{\mathcal{C}}

\newcommand{\GL}{{\rm GL}}
\newcommand{\SL}{{\rm SL}}
\newcommand{\Sps}{{\rm Sps}}

\renewcommand{\H}{\mathcal{H}}

\newcommand{\Aut}{{\rm Aut}}
\newcommand{\Int}{{\rm Int}}
\newcommand{\Hom}{{\rm Hom}}

\begin{document}

\title{On the almost algebraicity of groups of automorphisms of connected Lie groups}

\author{S.G.\ Dani}
\address{S.G.\ Dani\\
UM-DAE Centre for Excellence in Basic Sciences\\
University of Mumbai\\ 
 Vidyanagari, Santacruz\\ 
 Mumbai 400098, India}
 \email{shrigodani@cbs.ac.in, sdani.cebs@gmail.com}

 \author{Riddhi Shah}
\address{Riddhi Shah\\
School of Physical Sciences\\
Jawaharlal Nehru University\\ 
New Delhi 110067, India}
\email{rshah@jnu.ac.in, riddhi.kausti@gmail.com }

\begin{abstract}
Let $G$ be a connected Lie group, $C$ be the maximal compact connected  subgroup of the center of $G$, and let $\Aut (G)$ denote the group 
of Lie automorphisms of $G$, viewed, canonically, also as a subgroup of $\GL (\frak G)$, where $\frak G$ is the Lie algebra of $G$. 
It is known that when $C$ is trivial $\Aut (G)$ is almost algebraic, in the sense that it is of finite index in an algebraic subgroup 
of $\GL(\frak G)$, and in particular has only finitely many connected components. In this paper we analyse the situation further  in this respect, 
with $C$ possibly nontrivial, and describe necessary and sufficient conditions for almost algebraicity to hold; the criteria are  in terms of the group of 
restrictions of automorphisms of $G$ to $C$, and the abelian quotient Lie group $G/\overline{[G,G]}C$. For the class of Lie groups which 
admit a finite-dimensional representation with discrete kernel, a more specific criterion for $\Aut(G)$ to be almost algebraic is obtained, 
 while in the general case a variety of patterns are illustrated through examples. Along the way we also study almost algebraicity of subgroups 
 of $\Aut(G)$ fixing each point of a given torus in $G$, containing $C$, which also turns out to be of independent interest. 
\end{abstract} 

\maketitle



\tableofcontents

\section{Introduction}
Let $G$ be a connected Lie group, and $\frak G$ be the Lie algebra of $G$. We denote by $\Aut (G)$ the group of all Lie automorphisms of $G$, 
which we realize as a subgroup of $\GL ({\frak G})$, the general linear group of $\frak G$ over $\R$, by identifying each automorphism with its 
differential on $\frak G$. A subgroup of $\GL ({\frak G})$ is said to be an algebraic subgroup if it is of the form $\{\xi \in \GL ({\frak G})\mid p(\xi)=0\}$, 
where $p$ is a function on $\GL ({\frak G})$ which is a polynomial in the entries of $\xi$ and $(\det \xi)^{-1}$,   in a matrix representation with respect 
to a basis of $\frak G$.  A subgroup of $\GL ({\frak G})$ is said to be {\it almost algebraic} if it is a subgroup of finite index in an algebraic subgroup. 

It is of interest in various contexts to know  whether $\Aut (G)$, and various subgroups of it defined intrinsically, are algebraic or almost algebraic; 
see \cite{D-surv}. $\Aut (G)$ being almost algebraic implies in particular that it is almost connected, viz.\,it has only finitely connected components. On the other hand, by 
a theorem of Wigner \cite{W}  (see also \cite{D} and \cite{PW}) it is known that the connected component of the identity in $\Aut (G)$ is 
almost algebraic; thus  $\Aut (G)$ is almost algebraic  if and only if it is almost connected. Various conditions are known which 
ensure this; in particular it holds when the center of $G$ does not contain a compact subgroup of positive 
dimension (see \cite{D} and \cite{PW} for the latter result, and \cite{D-surv} for a broader perspective on the topic, and other related references). 

Here we  explore the theme further. Let $C$ denote the unique maximal compact connected subgroup of the center of $G$. We note that $C$ is a torus; let $d\geq 0$ be the dimension of $C$. Clearly $C$ is invariant under all $\alpha \in \Aut (G)$, and we get a natural homomorphism of $\rho:\Aut (G)\to \Aut (C)$, assigning to each automorphism of $G$ its restriction to $C$.  
We note that $\Aut (C)$ is isomorphic to $\GL(d, \Z)$, the group of $d\times d$ matrices with integer entries and determinant $\pm 1$;  it is a countable group, and it is  infinite when $d\geq 2$. For $\Aut (G)$ to be almost algebraic it has to be almost connected, and hence the image $\im \rho$ of $\rho$ has to be finite.  Thus when $d\geq 2$, in general $\Aut (G)$ need not be almost algebraic, as seen in particular when $G=\T^d$,  a torus of dimension~$d\geq 2$ or, more generally,  a connected abelian Lie group containing $\T^d$. $d\geq 2$, as a subgroup. 

{While the  condition that $\im \rho$ be finite is a necessary condition, it does not ensure almost algebraicity of $\Aut (G)$; see  the example of the `Walnut group'  in \S\,7.2. In  this paper we describe another necessary condition, and prove that the two conditions together are sufficient to ensure that $\Aut (G)$ is almost algebraic.

The other necessary condition  arises as follows.  Any (continuous)  homomorphism $\psi:G\to C$ whose restriction to $C$ is the identity map, induces an automorphism $\alpha$ of $G$ defined by $\alpha (g)=g\psi (g)$ for all $g\in G$. We call the automorphisms arising in this way {\it spiral shear automorphisms} (see \S\,5 for details); the term is inspired by the fact that analogously defined automorphisms with the center $Z$ of $G$ in place of $C$ have been called shear automorphisms (see \cite{D}, \cite{D-surv}). The spiral shear automorphisms form a closed subgroup of $\Aut (G)$ which we shall denote by  $\Sps (G)$.  We show that for $\Aut (G)$ to be almost connected, it is necessary that $\Sps (G)$  be  connected; cf.\  Theorem~\ref{prop-abel}.  It also turns out that the latter condition holds if and only if either $C$ is trivial or $G/\overline{[G,G]}C$ is simply connected; see Remark~\ref{rem1.2}. 

Our main result may now be stated as follows:

\begin{theorem}\label{characterization}
Let $G$ be a connected Lie group and $C$ be the maximal torus of the center of $G$. Then  $\Aut (G)$ is almost algebraic if and only if the following conditions are satisfied:
\begin{enumerate}
\item[{i)}] The subgroup of $\Aut (C)$ 
consisting of restrictions of automorphisms of $G$ to $C$ is a finite  group, and  

\item[{ii)}] Either $C$ is trivial or $G/\overline{[G,G]}C$ is simply connected. 
\end{enumerate}
\end{theorem} 

Now let $R$ be the radical of $G$ and $S$ be a semisimple Levi subgroup of $G$; (the statements made in the sequel in this regard do not depend on the choice of the semisimple Levi subgroup, they being conjugates of each other in $G$). With regard to the condition of $G/\overline{[G,G]}C$ being simply connected, as in~(ii) above, we note the following (proof~in~\S\,6). 

\begin{proposition}\label{rem:SR}
 Suppose that $S\cap R$ is finite. Then $G/\overline{[G,G]}C$ is simply connected if and only if $C$ is a maximal torus of $R$. 
\end{proposition}

Under an additional condition on the Lie group $G$  we deduce from the theorem a characterization for almost algebraicity, directly in terms of the structure of the Lie group (see Corollary~\ref{cor:charac}). A Lie group $G$ is said to be of 
{\it class} $\C$ if it admits a representation $\rho :G\to \GL(V)$, where $V$ is a finite-dimensional vector space over $\R$, such that 
$\ker \rho $ is discrete; the class includes in particular the linear Lie groups, as well as all simply connected Lie groups and all semisimple Lie 
groups. 
It is known (see \cite{DM}, Proposition~2.5) that $G$ is of class $\C$ if and only if $[R,R]$ is a closed simply connected (nilpotent) Lie subgroup, and also, if and only if $\overline {[R,R]}\cap Z$, where $Z$ is the 
center of $G$, has no compact subgroup of positive dimension.

\begin{corollary}\label{cor:charac}
Let $G$ and $C$ be as in Theorem~\ref{characterization} and suppose  that $G$ is of class~$\C$. 
Then we have the following:

1.  $\Aut (G)$ is almost algebraic if and only if either 
(i)~$C$~is trivial, or (ii)~$C$~is one-dimensional and  $G/\overline{[G,G]}C$ is simply connected.

2. If, moreover, $S\cap R$ is finite then $\Aut (G)$ is almost algebraic if and only if either 
(i)~$C$~is trivial, or (ii)~$C$~is one-dimensional and it is a maximal torus in $R$.   
\end{corollary}

When the Lie group $G$ is linear, viz.\,a subgroup of $\GL (d,\R)$ for some $d$, $S$ has finite center and hence $S\cap R$ is finite. Hence
Corollary~\ref{cor:charac}(2) generalizes the characterization of algebraicity of $\Aut (G)$ proved in \cite{CW}  for connected linear Lie groups $G$; in  \cite{CW} almost connected linear  Lie groups $G$ are also considered, obtaining results on almost algebaicity under additional  technical restrictions; we shall however not concern ourselves with that generality.

 \begin{remark}
 It may be noted that $\Aut (G)$ considered here, including in the case of linear Lie groups, consists of all continuous (and hence analytic) automorphisms of $G$. In the case of real algebraic groups (group of $\R$ points of an algebraic group defined over $\R$), which form a particular case of Lie groups, one may consider the group of their automorphisms that are (restrictions of) algebraic morphisms. Using the structure theory of these groups (cf.\  \cite{BT})  it can be readily  seen that the latter is almost connected if and only if the radical of $G$ coincides with the unipotent radical. We shall not go into the details here. 
 \end{remark} 

In the course of proof of Theorem~\ref{characterization} we also consider the question of almost algebraicity  of subgroups of $\Aut(G)$ which fix 
pointwise a given torus (compact connected abelian subgroup) in $G$. Consideration of this issue was motivated by the study in \cite{CS}, dealing 
with a characterization of automorphisms of a connected Lie group  acting distally on certain  compact spaces.  On the other hand it also plays 
a crucial role in the proof of Theorem~\ref{characterization}. 

For a torus $T$ in $G$  we denote by $F_T(G)$ the 
subgroup of $\Aut (G)$ consisting of all  $\tau \in \Aut(G)$ such that $\tau (x)=x$ for all $x\in T$. We prove the following. 

\begin{theorem}\label{torus-fixing}
Let $G$ and $C$ be as in Theorem \ref{characterization}, and let $T$ be a torus in $G$ 
containing $C$.  Then $F_T(G)$ is almost algebraic if and only if either $C$ is trivial or $G/ \overline {[G,G]}T$  is simply connected. 
In particular, if $T$ is a maximal torus in $G$ then $F_T(G)$ is almost algebraic. Also, if $S\cap R$ is finite and $T$ is a maximal torus in $R$ then $F_T(G)$ is almost algebraic.
\end{theorem} 

We note that for a torus $T$ {\it not} containing $C$, $F_T(G)$ need not be almost algebraic even if $G/ \overline {[G,G]}T$  is simply connected;  see \S\,7.1.3, for an example. 

The paper is organized as follows. In \S\,2 we introduce notations and discuss some preliminaries. \S\,3 is devoted to the proof of the sufficiency parts of Theorems~\ref{characterization} and~\ref{torus-fixing}.  In \S\,4 
we present some general results on the structure of $\Aut (G)$, that are reminiscent of some results 
from \cite{D}, together with some strengthenings required here.   In  \S\,5 we discuss the  class of ``spiral shear automorphisms"; these are the source of  the other obstruction to 
almost algebraicity of $\Aut (G)$  involved in Theorem~\ref{characterization}. Proofs of the results stated above are completed 
 in \S\,6. In the last section, \S\,7, we discuss a wide range of examples
illustrating in particular the relevance of the conditions, and interdependence  of various conditions involved, in the results described above.

\section{Notations and preliminaries}

The notation $\Aut (.)$ as above will be used, in place of $G$ as above, for any connected Lie group, or connected Lie subgroup $H$; thus
$\Aut (H)$ will denote the group of all Lie automorphisms of $H$, realized as a subgroup of $\GL ({\frak H})$, with $\frak H$ the Lie algebra of $H$. In terms of the 
theory of algebraic groups, the group  $\GL ({\frak H})$ is also considered as the group of $\R$-elements of $\bf GL$, the general linear algebraic group 
corresponding to the dimension of $\frak H$.  A subgroup of $\GL ({\frak H})$ is thus an algebraic subgroup if it is the group of $\R$-points of an algebraic 
subgroup of $\bf GL$; this notion coincides with that of algebraic subgroups of $\Aut (G)$ referred to in the introduction. 

Now let $G$ be a connected Lie group and $\frak G$ denote the Lie algebra of $G$. 
Recall that $\Aut (G)$ is in particular a Lie group which we view, canonically, as a subgroup of $\GL (\frak G)$, identifying each $\alpha \in \Aut (G)$ with $d\alpha \in \GL (\frak G)$. 
For any Lie group or Lie subgroup $S$ we shall denote by $S^0$ the connected component 
of the identity element in $S$. 

We denote by 
$\tilde G$ a (simply connected) univeral covering  group of $G$ and by $\eta : \tilde G \to G$ the covering homomorphism. The Lie algebra of $\tilde G$ will be identified with $\frak G$, canonically so that $d\eta $ is the identity automorphism of $\frak G$.  
We denote by $\Aut (\frak G)$ the subgroup of $\GL(\frak G)$ consisting of Lie algebra automorphisms of $\frak G$. For a connected Lie subgroup $H$ of $G$ we denote by $\tilde H$ the connected Lie subgroup $\eta^{-1}(H)^0$. 

For any 
closed subgroup $H$ of $G$ we shall denote by $\Aut_H(G)$ the subgroup consisting of all $\alpha \in \Aut (G)$ such that  
$\alpha (H)=H$.
Also, for any 
closed subgroup $H$ of $G$ we shall denote 
by $F_H(G)$ the subgroup consisting of all $\alpha \in \Aut (G)$ such that  $\alpha (h)=h$ for all $h\in H$. 

 For any $x\in G$ we denote by $\sigma_x:G\to G$ the inner automorphism associated with $x$, defined by $\sigma_x(g)=xgx^{-1}$ for all 
 $g\in G$.  Also, for any subgroup $H$ of $G$ we shall 
 denote by $\Int_H(G)$ the subgroup $\{\sigma_s\mid s\in H\}$ of $\Aut (G)$ (viewed as a subgroup of $\GL (\frak G)$). We denote by $Z$ the center of $G$. We note that if for a subgroup $H$ of $G$, $\Int_H(G)$ is a closed subgroup of $\GL (\frak G)$ then $HZ$ is a closed subgroup of $G$.

The following remark plays a crucial role in many arguments in the sequel.

\begin{remark}\label{rem-sc}
  Under the identifications as above $\Aut (G)\subseteq \Aut (\frak G)$, and equality holds when $G$ is simply connected, since in this case every Lie algebra automorphism of $\frak G$ corresponds to a Lie automorphism of $G$.  Thus we have $\Aut (G)\subseteq \Aut (\frak G)=\Aut (\tilde G)$. Since the conditions  for a linear transformation of $\frak G$ to be a Lie algebra automorphism of $\frak G$  are polynomial identities, $\Aut (\frak G)$  is an algebraic subgroup of $\GL(\frak G)$;   in particular,  $\Aut (\tilde G)$ is an algebraic subgroup of  $\GL(\frak G)$. Moreover, for any connected Lie subgroup $H$ of $\tilde G$, $F_H(\tilde G)$ is an  algebraic subgroup, as it consists of all $\alpha \in \Aut (\tilde G)$ such that every element of the Lie subalgebra $\frak S$ of $S$ is fixed by $d\alpha$.  
  Secondly, as a subgroup of $\Aut (\tilde G)$, $\Aut (G)$ consists precisely of $\{\tau \in \Aut (\tilde G)\mid \tau (\ker \eta)=\ker \eta\}$. Thus the task of showing that $\Aut (G)$ is almost algebraic is equivalent to showing that the latter subgroup is an almost algebraic subgroup of $\Aut (\tilde G)$. 
\end{remark}

We denote by $R$ the radical of $G$, by $N$ the nilradical of $G$,  and $S$ a semisimple Levi subgroup of $G$. We note that $R$ and $N$ are closed normal subgroups of $G$, while $S$ is a connected Lie subgroup of $G$ which may in general be neither normal nor closed. Since $S$ and $N$ are connected Lie subgroups of $G$ and $N$ is a normal subgroup, it follows that $SN$ is a Lie subgroup of $G$; it turns out that $SN$ is a closed subgroup (see Proposition ~\ref{SN}). Also, since $G/N$ is reductive, $SN/N$ is normal in $G/N$ it follows that $SN$ is a normal subgroup of $G$.

\section{Sufficiency of the conditions for algebraicity of $\Aut (G)$}

We  now suppose that $G$ and $C$ are as  in Theorem~\ref{characterization}, and   conditions (i) and (ii) are satisfied; the proof of almost algebraicity of $G$  will be completed in this section (see Corollary~\ref{direct}), after some preparatory steps. We begin by noting the following.

 \begin{proposition}\label{SN}
 For any  connected Lie subgroup $M$ of $N$ which contains $Z^0$ and is normalized by $S$, $SM$ is a closed subgroup. 
 In particular  $SN$ is a closed normal subgroup containing $\overline{[G,G]}C$. 
 \end{proposition}
\proof  
As $S$ is semisimple, $\Int_S(G)$ is an almost algebraic subgroup of $\GL (\frak G)$. Now let $M$ be as in the statement. Then  $\Int_M(G)$ is an almost algebraic subgroup, as it is a connected Lie subgroup consisting of unipotent elements. Therefore $\Int_{SM}(G)$ which is the same as $\Int_S(G) \Int_M(G)$ is almost algebraic, and in particular it is a closed subgroup of $\GL (\frak G)$. This implies that $SMZ$ is a closed subgroup of $G$. 
Since $Z^0\subset M$, $SM$ is the connected component of the identity in $SMZ$, and hence it is closed.  
In particular $SN$ is closed. As it contains $[G,G]$ it is a normal subgroup, and, $SN$ being closed,  we get  $\overline{[G,G]}\subset SN$. Since $C\subset N$ we now get that $\overline{[G,G]}C \subset SN$. \qed

\begin{proposition}\label{center}
Let $M$ be a connected Lie subgroup of $N$ containing $C$ and normalized by $S$. Let $T_0$ be a torus in $G$ such that $T_0$  centralizes $S$, normalizes $M$, and $T_0\cap SM$ is finite. Let $Q=T_0S$.
 Let  $Z'$ be the center of $\tilde Q\tilde M$. If $q\in \tilde Q$ and $m\in \tilde M$ are such that   $qm\in Z'$, then $q$ and $ m$ are contained in $Z'$. Moreover,  $Z'\cap \tilde Q$ is fixed pointwise  by all $\alpha \in F_{\tilde T_0}(\tilde Q\tilde M)^0$.
\end{proposition}
\proof We note that under the conditions in the statement $\tilde Q\tilde M$ is a semidirect product of $\tilde Q$ and $\tilde M$. We
prove the first assertion  by induction on the dimension of $M$; the statement is  obvious when the dimension is $0$, or $1$, and we may assume it to be true in all cases when $M$ is of lower dimension than under consideration.  If the action of $\tilde Q$ on $\tilde M$ is trivial, then $\tilde Q \tilde M$ is a direct product of $\tilde Q$ and $\tilde M$, in which case the desired assertion holds. Suppose the action is nontrivial. Since $\tilde M$ is a nilpotent Lie group this implies that the  $\tilde Q$-action on $\tilde M/[\tilde M,\tilde M]$ is nontrivial.  It follows that  there exists a $\tilde Q$-invariant closed connected Lie subgroup $\tilde M'$ of $\tilde M$ containing $[\tilde M,\tilde M]$, such that  the quotient $\tilde M/\tilde M'$ is a vector space, and  the $\tilde Q$-action on it is nontrivial and irreducible. It follows therefore that  $Z'$ is contained in $\tilde Q\tilde M'$. The desired assertion then follows from the induction hypothesis. 

Now let $D$ be the subgroup consisting of all $x\in \tilde S$ such that $x\tilde T_0\tilde M \cap Z'$ is nonempty. Since $\tilde T_0\tilde M$ and $ Z'$ are invariant the action of $\Aut (\tilde Q\tilde M)$, we get that $D$ is also $\Aut (\tilde Q\tilde M)$-invariant. On the other hand $D$ is contained in the center of $\tilde S$ and hence it is discrete. It follows that $D$ is pointwise fixed by any automorphism $\alpha \in \Aut (\tilde Q\tilde M)^0$. We note that $Z'\cap \tilde Q$ is contained in $D\tilde T_0$. Therefore all elements of $Z'\cap \tilde Q$ are fixed by any $\alpha \in F_{\tilde T_0}(\tilde Q\tilde M)^0$. 
 \qed

The following Corolary, together with the fact that $\Aut (G)$ is almost algebraic when $C$ is trivial (cf.\,\cite{D} and \cite{PW}),  establishes the ``if'' parts of Theorems~\ref{characterization} and \ref{torus-fixing}.

\begin{corollary}\label{direct}
Let $G$ and $C$ be as in Theorem~\ref{characterization}. Let $T$ be a torus in $G$ containing $C$ such  that $ G/\overline {[G,G]}T$ is simply connected. Then  $F_T(G)$ is  almost algebraic.  Moreover, if   $ G/\overline {[G,G]}C$ is simply connected and Condition~(i) as in Theorem~\ref{characterization} holds then $\Aut (G)$ is almost algebraic. 
\end{corollary} 

\proof We recall that $SN$ is a closed normal subgroup (see Proposition~\ref{SN}), and set $H=TSN$.   Considering the conjugation action of $S$ on $H$ we see that there exists a torus $T_0$ contained in $T$ and centralized by $S$ and such that $T_0\cap SN$ is finite and $T_0SN=TSN=H$. 

Now  let $Q=T_0S$;  thus $H=QN$. By Proposition~\ref{SN} $\overline{[G,G]}C$ is contained in $SN$ and so $\overline{[G,G]}T$ is contained in $H$.  As $H$ is connected and by  hypothesis $ G/\overline {[G,G]}T$ is simply connected, this implies  that $G/H$ is simply connected.   
Since 
$H$ is a closed normal subgroup of $G$ this implies that  $\ker \eta$ is contained in $\tilde H$. Let $Z'$ be the center of $\tilde H$. Then $\ker \eta \subset Z'$.   Let $\tilde M$ be the subgroup of $\tilde N$ consisting of all $x\in \tilde N$ which are centralized by $\tilde  Q$. Then $\tilde M$  and $\tilde Q\tilde M$ are closed connected subgroups. We note that $Z'\subset \tilde Q\tilde M$ and in particular $\ker \eta \subset \tilde Q\tilde M$. Now $\tilde Q\tilde C$ is a closed normal subgroup of $\tilde Q\tilde M$ and $\tilde Q\tilde M/\tilde Q\tilde C$ is simply connected. As $\ker \eta \subset \tilde Q\tilde M$ this further implies that $\ker \eta \subset \tilde Q\tilde C$. Also, since $\ker \eta \subset Z'$ and $\tilde C\subset Z'$, we get that $\ker \eta \subset (Z'\cap \tilde Q)\tilde C$. 

By Proposition~\ref{center}, applied choosing $M=N$,  we have that $Z'\cap \tilde Q$ is pointwise fixed by all $\alpha \in F_{\tilde T_0}(\tilde H)^0$. The subgroup $\tilde S\tilde N$ is invariant under all automorphisms of $\tilde G$ and hence the subgroup $\tilde H$ which equals $\tilde T_0\tilde S\tilde N$ is invariant under any 
$\alpha \in F_{\tilde T_0}(\tilde G)$.  Hence the preceding conclusion implies  that $Z'\cap \tilde Q$ is pointwise fixed by all $\alpha \in F_{\tilde T_0}(\tilde G)^0$. Since $T_0$ and $C$ are contained in $T$, and $\ker \eta \subset (Z'\cap \tilde Q)\tilde C$ we get that $\ker \eta $ is pointwise fixed by all  $\alpha \in F_{\tilde T}(\tilde G)^0$.
 Hence,  in  $\GL (\frak G)$ we have, $F_{\tilde T}(\tilde G)^0\subset F_T(G)\subset F_{\tilde T}(\tilde G)$. 
As $F_{\tilde T}(\tilde G)$ is almost algebraic this implies that $F_T(G)$ is almost algebraic.  This proves the first assertion of the Corollary. 

Now suppose that  $ G/\overline {[G,G]}C$ is simply connected. Then by the above conclusion $F_C(G)$ is almost algebraic. 
Since $F_C(G)$ is the kernel of the homomorphism $\rho:\Aut (G) \to\Aut (C)$,  assigning to each automorphism its restriction to $C$, it follows that when Condition~(i) of Theorem~\ref{characterization} holds $\Aut (G)$ is almost algebraic. \qed

\section{Structure of the automorphism group}

In this section we prove some results on the structure of the automorphism groups of Lie groups which will be used in the sequel. We follow the notations introduced earlier. 
 We begin by recalling  the following useful fact. 

\begin{theorem}\label{Chevalley}
Let  $H$ be any closed connected subgroup of $G$. Then $\Int_{[H,H]}(G)$ is an almost algebraic subgroup. 
In particular it is a closed subgroup, and $\Int_{\overline{ [H,H]}}(G)=\Int_{[H,H]}(G)$.
\end{theorem}

\begin{proof} The first statement is immediate from  \cite{C}, Ch.\,II, Theorem 15. It is clear that subgroups that are almost algebraic,  are closed subgroups, 
and therefore the last statement  in the theorem holds. \end{proof}

In particular $\Int_{[G,G]}(G)$ is an almost algebraic group; it may be noted that  in general $\Int_G(G)$ need not be almost algebraic.  For example, if $G$ is the semidirect product of the one-parameter subgroup $\{{\rm diag}\,(e^{\lambda_1t}, \dots , e^{\lambda_nt})\mid t\in \R\}$ with $\R^n$, $n\geq 2$, then $\Int_G( G)$ is almost algebraic if and only if $(\lambda_1, \dots, \lambda_n)$ is a scalar multiple of an integral vector. 

The following result is a generalization of Lemma 3 in \cite{D}, together with a strengthening (see Remark~\ref{rem2.3} below). 

\begin{lemma}\label{lem1} 
Let  $H$ be a closed connected subgroup of $G$ 
and $T$ be a maximal torus in $H$. Then $$\Aut_H (G)=\Int_{{[H.H]}}(G)\cdot (\Aut_T(G)\cap \Aut_H(G)).$$ 
\end{lemma} 

\begin{proof} Consider $\alpha \in \Aut_H(G)$. Then  $T$ and $\alpha (T)$ are maximal tori in $\overline {[H.H]}T$. As the latter is a connected 
Lie group, the two subgroups are conjugate in $\overline {[H.H]}T$; thus  there exists $h\in \overline {[H.H]}$ such that 
$\alpha (T)=hTh^{-1}=\sigma_h (T)$. Let $\tau = \sigma_h^{-1}\circ \alpha$. Then $\tau \in \Aut_T(G)\cap \Aut_H(G)$ and $\alpha =\sigma_h\circ \tau$.
As $\sigma_h\in \Int_{\overline{ [H,H]}}(G)=\Int_{[H,H]}(G)$, 
by Theorem~\ref{Chevalley}, this proves the Lemma. \end{proof}

\begin{remark}\label{rem2.3}
When $H=G$ the Lemma reduces to $\Aut (G)=\Int_{{[G.G]}}(G)\cdot \Aut_T(G)$. Such a decomposition was noted in 
\cite{D} with respect to a chosen subgroup containing $[G,G]$, denoted there by $H$, in place of ${[G,G]}$ in the above statement;  see also a 
clarification/correction given in \cite{PW} towards ascertaining the validity of the statement.  As ${[G,G]}\subset H$ the present assertion is stronger than Lemma~3 in \cite{D}, even in the special case with $H=G$.  
\end{remark}

\begin{corollary}\label{cor2.4}
Let  $H$ be a closed connected subgroup of $G$ and $T$ be a maximal torus in $H$. Let $\H$ be a subgroup 
of  $\Aut_H(G)$ containing $\Int_{[H,H]} (G)$. Then $\H$ is almost algebraic if and only if $\Aut_T(G)\cap \H$ is almost algebraic. 
\end{corollary}

\begin{proof}   
By Lemma~\ref{lem1} ${\H} =\Int_{[H,H]}(G)\cdot (\Aut_T(G)\cap {\H})$ and by  Theorem~\ref{Chevalley}
$\Int_{[H,H]}(G)$ is almost algebraic. Therefore it follows that if $\Aut_T(G)\cap \H$ is almost algebraic then $\H$ is almost algebraic. 
Now suppose that $\H$ is almost algebraic. We note that $\Aut_T(G)\cap \H$ consists precisely of all $\tau \in \H$ such that 
$\tau (T)=T$, and the latter condition is equivalent to $d\tau ({\frak T}) ={\frak T}$, where $\frak T$ is the Lie subalgebra 
corresponding to $T$ and $d\tau$ is the differential of $\tau$. This shows that $\Aut_T(G)\cap \H$ is almost algebraic. \end{proof} 

\begin{corollary}\label{cor2.5}
Let $T$ be a maximal torus in $G$ and let $N_{[G,G]}(T)$ and $Z_{[G,G]}(T)$ denote respectively the normalizer and centralizer of $T$ in $[G,G]$.  Then $N_{[G,G]}(T)/Z_{[G,G]}(T)$ is finite. In particular if $g\in [G,G]$ is such that $T$ is $\sigma_g$-invariant, then every orbit of the $\sigma_g$-action on $T$ is finite.

\end{corollary}

\begin{proof}   
For any $g\in N_{[G,G]}(T)$ we get canonically an automorphism from  $\Int_{[G,G]}(G)\cap \Aut_T(G)$ and the automorphism is trivial if and only if $g\in Z_{[G,G]}(T)$. Since $\Aut (T)$ is countable the action of any $\tau \in (\Int_{[G,G]}(G)\cap \Aut_T(G))^0$ on $T$ is trivial.  
Thus the group $N_{[G,G]}(T)/Z_{[G,G]}(T)$ is isomorphic to a quotient of $(\Int_{[G,G]}(G)\cap \Aut_T(G))/(\Int_{[G,G]}(G)\cap \Aut_T(G))^0$. 
From Corollary~\ref{cor2.4}, choosing $H=G$ and ${\H}=\Int_{[G,G]}(G)$  we get that $\Int_{[G,G]}(G)\cap \Aut_T(G)$ 
is almost algebraic, and hence $(\Int_{[G,G]}(G)\cap \Aut_T(G))/(\Int_{[G,G]}(G)\cap \Aut_T(G))^0$ is finite. Therefore $N_{[G,G]}(T)/Z_{[G,G]}(T)$ is finite.
The second assertion follows immediately from the first.  
\end{proof}

We note that in general $N_G(T)/Z_G(T)$ need not be finite. Consider the semidirect product $G=P\cdot T$ where $P$ is a one-parameter subgroup $\{P(t)\mid t\in \R\}$, $T$ is the torus $\{(z_1,z_2)\in \Cm^2\mid |z_1|=|z_2|=1\}$, with the action of $P$ on $T$  given by  $P(t)(z_1,z_2)=(z_1, e^{i\alpha t}z_2)$ for all $t\in \R$ and $z_1, z_2\in T$, where $\alpha$ is an irrational number.   
It is readily seen that, $N_G(T)=G$ and $Z_G(T)=T$, so $N_G(T)/Z_T(G)$ is infinite, in fact, isomorphic to $\R$. 

\begin{corollary}\label{cor:finite}
Let $\tau \in \Aut (G)^0$ and $T$ be maximal torus invariant under the action of $\tau$. Then the orbits of the $\tau$-action on $T$ are finite. 
\end{corollary}

\proof 
By Lemma~\ref{lem1} we have $\Aut (G)=\Int_{{[G.G]}}(G)\cdot \Aut_T(G)$, and hence $\Aut (G)^0=\Int_{{[G.G]}}(G)\cdot \Aut_T(G)^0$. 
Hence there exist $g\in [G,G]$ and $\alpha \in \Aut_T(G)^0$ such that $\tau =\sigma_g\circ \alpha$.   As $T$ is invariant under $\tau$ and $\alpha$, it is also invariant under $\sigma_g$. Moreover, since $\alpha \in \Aut_T (G)^0$, its action on $T$ is trivial. Therefore  the orbits of $\tau$ and $\sigma_g$ on $T$ coincide, and hence by Corolllary~\ref{cor2.5} they are finite. \qed

\section{Spiral shear automorphisms}

As before, $G$ is a connected Lie group and $C$  the maximal torus of  the center of $G$.  
Recall that $\Sps (G)$ denotes the group of all spiral shear automorphisms of $G$; (see \S\,1). Let  $A=G/\overline{[G,G]}C$, which is a connected abelian Lie group. Let  $\Hom (A,C)$ denote the group of all continuous homomorphisms from $A$ to $C$. For each $\varphi \in \Hom(A,C)$ we have a spiral shear automorphism given by $\Theta (\varphi)(g)=g\varphi (g\overline{[G,G]}C)$ for all $g\in G$, which  we shall denote by $\Theta (\varphi)$. We note also that if $\psi : G\to C$ is a continuous homomorhphism whose restriction to $C$ is the identity map (in the way the transformations were introduced in the introduction) then $\psi $ is trivial on 
$\overline {[G,G]}C$ and factors to a homomorphism $\varphi: A\to C$ as above. 

We equip the group $\Hom(A,C)$ with the compact open topology (corresponding to unifrom convergence on compact subsets).  Let $\Theta :\Hom (A,C)\to \Aut (G)$ be the map defined by  $\varphi \mapsto \Theta (\varphi)$ as above; its  range is $\Sps (G)$.  It is straightforward to verify that $\Theta$ is a continuous injective homomorphism of the topological groups. Hence we get a  homomorphism $\bar \Theta: \Hom (A,C)/\Hom (A,C)^0 \to \Aut (G)/\Aut (G)^0$, induced by $\Theta$.

\begin{remark}\label{rem1.2}
As  a connected abelian Lie group $A$ is  a direct product $\R^n \times \T^m$, $n,m\geq 0$, where $\T=\R/\Z$ is the circle group.   We see that $\Hom (A,C)$ is topologically isomorphic to 
$\R^{nc}\times \Z^{mc}$, where $c$ is the dimension of $C$. In particular, Hom$(A,C)$ is connected if and only if $mc=0$, namely 
if either $C$  is trivial or $A$  is $\R^n$. We also note  that $\Hom (A,C)/\Hom (A,C)^0$ is isomorphic to $\Z^{mc}$, and if $mc>0$ it has no nontrivial elements of finite order. 
\end{remark}

 The following theorem is a crucial ingredient towards proving the ``only if'' part of  Theorem\ref{characterization}.

\begin{theorem}\label{prop-abel}
Let $K$ be the maximal compact subgroup of $A$. Let  $\varphi \in {\rm Hom} (A,C)$.  Then $\Theta (\varphi)$ is contained in $\Aut (G)^0$  if and only if $\varphi (K)$  is the trivial subgroup.  Consequently,

(i) $\bar \Theta:  {\rm Hom}(A,C)/ {\rm Hom}(A,C)^0 \to \Aut (G)/\Aut (G)^0$ is injective, and 

(ii) if $C$ and $K$ are nontrivial then there exists 
$\alpha \in \Sps (G)$ such that   $\alpha^n$ is not contained in $\Aut (G)^0$ for any $n\in \Z$. 
\end{theorem} 

\begin{proof}   From the structure theory of connected 
abelian Lie groups we get that  $\Hom (A,C)^0$ consists precisely of $\varphi \in {\rm Hom}\,(A,C)$ such that 
$\varphi  (K)$ is trivial. Since $\Theta$ is continuous, this implies that if $\varphi  (K)$ is trivial then $\Theta (\varphi )$ is contained 
in $\Aut (G)^0$. 

Now consider any $\varphi \in {\rm Hom} (A,C)$ such that  $\varphi (K)$ is nontrivial. Then  $\varphi (K)$ is a torus of positive dimension and there exists $k\in K$ such that $\varphi (k)$ is an element of infinite order in $C$.  Let $T$ be a maximal torus in $G$. As it contains $C$,
 $T$ is invariant under the action of $\Theta (\varphi)$, and the  $\Theta (\varphi)$-orbit of $t\in T$ consists of $\{t\varphi (k)^j\mid j\in \Z\}$. Since $\varphi (k)$ is of infinite order, the orbits are infinite. By Corollary~\ref{cor:finite} this implies that $\Theta (\varphi)$ is not contained in $\Aut (G)^0$.   This proves the initial assertion in the theorem. 

As Hom$(A,C)^0$ consists precisely of all $\varphi$ such that $\varphi (K)$ is trivial,  the conclusion above means that for $\varphi\in {\rm Hom}(A,C)$, $\Theta (\varphi)\in  \Aut (G)^0$ if and only if $\varphi \in  {\rm Hom}(A,C)^0$. 
This proves~(i). Now suppose that $C$ and $K$ are nontrivial. Then there exists $\varphi \in {\rm Hom}(A,C)$ such that $\varphi^n\notin {\rm Hom}(A,C)^0$ for any  $n\in \Z$ (cf.\ Remark~\ref{rem1.2}). Hence for $\alpha = \Theta (\varphi)\in \Sps (G)$,  $\alpha^n$ is not contained in $\Aut (G)^0$ for any $n\in \Z$. This proves (ii).\end{proof}

\section{Proofs of the main results}

\begin{proof}[Proof of Theorem~\ref{characterization}] As noted before, Corollary~\ref{direct} establishes the ``if'' part. We now prove the converse.  Suppose that $\Aut (G)$ is almost algebraic.  Then, as noted in the introduction,  Condition~(i) must hold since for the
homomorphism of $\rho:\Aut (G)\to \Aut (C)$, assigning to each automorphism of $G$ its restriction to $C$,  $\im \rho$ has to be almost connected, and since $\Aut (C)$ is countable it has to be finite. To complete the proof it suffices to prove that if $C$ is nontrivial then $G/\overline{[G,G]C}$ is simply connected. As $\Aut (G)$ is almost connected,  
for any $\alpha \in \Aut (G)$ there exists $n\in \N$ such that $\alpha^n\in \Aut(G)^0$. As $C$ is assumed to be nontrivial, by 
Theorem~\ref{prop-abel} this implies that  $G/\overline{[G,G]}C$ does not contain a torus of positive dimension; hence it is simply 
connected.  \end{proof}

\begin{proof}[Proof of Proposition~\ref{rem:SR}]
Let $T$ be a maximal torus in $G$. As $S\cap R$ is finite, $T$ can be expressed as 
 $T=(T\cap S)(T\cap R)$. Thus if $C$ is a maximal torus in $R$ we get that $SC$ contains a maximal torus in $G$ and since $SC$ is contained in 
$\overline{[G,G]}C$ and the latter is a closed normal subgroup of $G$, it follows that 
$G/\overline{[G,G]}C$ is simply connected. 

Conversely suppose that $G/\overline{[G,G]}C$ is simply connected. Recall that $SN$, where $N$ is the nilradical, is a closed normal subgroup containing $\overline{[G,G]}C$ (cf.  Proposition~\ref{SN}), and hence we get that $G/SN$ is simply connected. Let $T$ be a torus in $R$. As $G/SN$ is simply connected, we then get that $T\subset SN$ and since $T\subset R$ this means that $T\subset N$. Since $C$ is the unique maximal torus in $N$ it follows that $T$ is contained in $C$. Thus $C$ is a maximal torus in $R$.\end{proof}

\medskip
Towards our proof of Corollary~\ref{cor:charac} we first note the following:

\begin{proposition}\label{classC}
If $G$ is of class $\C$, then $\overline{[G,G]}\cap R$ contains no nontrivial torus.
\end{proposition}

\begin{proof}  It is noted in Proposition~2.5 of \cite{DM} that when $G$ is of class $\C$, $[R,R]$ 
is a closed simply connected nilpotent Lie subgroup, and in particular contains no nontrivial compact subgroup; (the proof is by direct application of 
Lie's theorem). It therefore suffices to prove the statement for $G/[R,R]$ (with $R/[R,R]$ as its radical), in place of $G$ and $R$. We may therefore 
assume $R$ to be abelian.  Then $R$ may be decomposed as $V\times T$, where $V$ is a vector 
group, $T$ is a torus, and $V$ and $T$ are invariant under the conjugation action of the semisimple Levi subgroup $S$. We note that, as the automorphism group of 
$T$ is countable, the action of $S$ on $T$ is trivial. It follows that $\overline{[G,G]}=SW$, where $W$  is the vector subspace of $V$ consisting 
of the sum of all irreducible $S$-invariant subspaces of $V$ on which the action is nontrivial. Thus $\overline{[G,G]}\cap R= (S\cap R)W$. 
This shows that it contains no nontrivial torus. \end{proof}

\medskip

\begin{proposition}\label{lem5.5}
Suppose that 
$\overline{[G,G]}\cap C$ is finite. Then there exists a closed normal subgroup $G_1$ of $G$ such that $G$ is an almost direct product of $G_1$ 
and $C$. 
The conclusion holds in particular if $G$ is of class~$\C$. 
\end{proposition}

\begin{proof} Let $A=G/\overline{[G,G]}$ and $\pi:G\to A$ be the canonical quotient homomorphism.  Then $\pi (C)$ is a torus in $A$. As $A$ is 
an abelian Lie group there exists a closed subgroup $A_1$ of $A$ such that $A$ is a direct product of $A_1$ and $\pi (C)$. Let 
$G_1=\pi^{-1}(A_1)$. Then $G=G_1C$. As $G_1\cap C=\overline{[G,G]}\cap C$ is finite, it follows that $G$ is 
an almost direct product of $G_1$ and $C$. The second assertion now follows from 
Proposition~\ref{classC}. \end{proof}

\medskip
\begin{proof}[Proof of Corollary~\ref{cor:charac}] We follow the notation as in the statement. As recalled before, if $C$ is trivial then $\Aut (G)$ 
is almost algebraic (cf.\,\cite{D} and \cite{PW}). When  $C$ is one-dimensional $\Aut (C)$ is finite, and under the hypothesis that  
$G/\overline{[G,G]}C$ is simply connected, Theorem~\ref{characterization} implies that $\Aut (G)$ is almost algebraic. 
Also, if $C$ is nontrivial and $G/\overline {[G,G]}C$ contains a nontrivial torus, then by Theorem~\ref{prop-abel}  $\Aut (G)$ is not almost algebraic. 
(It may be noted that the proof of the part so far does not involve the condition that $G$ is of class~$\C$). To prove assertion~(1) in the statement it remains to show that  when $G$ is of 
class~$\C$, and $C$ is of dimension greater than 1 then $\Aut (G)$ is not almost algebraic. 

Suppose  that $G$ is of class $\C$ and $C$ is of dimension greater than $1$.  
By Lemma~\ref{lem5.5} there exists a closed connected normal subgroup $G_1$ of $G$ such that $G$ is an almost direct 
product of $G_1$ and $C$.  For any $\tau \in \Aut (C)$ which fixes $G_1\cap C$ pointwise, we have a well-defined automorphism $\alpha \in \Aut (G_1\times C)$ given by $\alpha (xc)=x\tau (c)$ for all $x\in G_1$ and $c\in C$. This shows that $\im \rho$ (where $\rho:\Aut (G)\to \Aut (C)$ is the restriction map as before) contains all $\tau \in \Aut (C)$ which fix $G_1\cap C$ pointwise. The latter is a subgroup of finite index in $\Aut (C)$ and hence it is infinite, as $C$ is of dimension greater than $1$. This shows that $\Aut (G)$ is not almost algebraic, thus proving~(1). Assertion~(2) follows directly from assertion~(1) and Proposition~\ref{rem:SR}.
\end{proof}

 \begin{proof}[Proof of  Theorem~\ref{torus-fixing}] 
Let $G$, $C$ and $T$ be as in the statement. If $C$ is trivial then $\Aut (G) $ is almost algebraic (cf.\,\cite{D} and \cite{PW}), and hence $F_T(G)$ is almost algebraic; see Remark~\ref{rem-sc}. In the case when $G/\overline {[G,G]}T$ is simply connected Corollary~\ref{direct} implies that $F_T(G)$ is almost algebraic. This proves the ``if'' part of the theorem. 
	 
	 Now suppose that  $C$ is nontrivial and  $G/\overline{[G,G]}T$ is not simply connected. As in \S\,5 let $A=G/\overline{[G,G]}C$ and let  $K$ be the maximal compact subgroup of $A$. Let $B =\overline{[G,G]}T/\overline{[G,G]}C$. As $G/\overline{[G,G]}T$ is not simply connected, $K$ is not contained in $B$. Since $C$ is nontrivial we can get a $\varphi \in \Hom (A,C)$ such that $\varphi (B)$ is trivial and $\varphi (K)$ is nontrivial. Let $\alpha =\Theta (\varphi)$, in the notation of \S\,5. Since $\varphi (B)$ is trivial it follows that $\alpha \in F_T(G) $.  On the other hand, since $\varphi (K)$ is nontrivial,  by Theorem~\ref{prop-abel} 
the image of $\alpha$ in the quotient	 $\Aut (G)/\Aut (G)^0$ is of infinite order. This shows that $F_T(G)$ is not almost connected and hence not almost algebraic, which proves the only if part.  

We note that if either $T$ is a maximal torus in $G$, or $S\cap R$ is finite and $T$ is a maximal torus in $R$, then $G/\overline{[G,G]}T$ is simply connected, which implies the remaining assertions. 
 \end{proof}

\section{Examples and concluding remarks}

In this section we discuss some special classes of Lie groups, and exhibit examples   bringing out  the significance of various conditions involved in the results described.  

\subsection{Automorphism groups of nilpotent Lie groups}

 To begin with consider  any connected nilpotent Lie group $N$. Let $C$ be the maximal compact subgroup of the center 
of $N$. Then  $C$ is the maximal torus in $N$. Hence in this case by Corollary~\ref{direct} $F_C(N)$ is almost algebraic. 

Concerning algebraicity of $\Aut (N)$ however there are several possibilities. Let $d$ be the dimension of $C$. If $d\leq 1$ then by Theorem~\ref{characterization} $\Aut(G)$ is 
almost algebraic, since in this case $\Aut (C)$ is finite and $N/C$ is simply connected. Now suppose $d\geq 2$. 
If $\overline{[N,N]}\cap C$ is trivial, which holds in particular if $N$ is abelian, then by Proposition~\ref{lem5.5}  
$N$  is a direct product of a simply connected nilpotent Lie group $N_1$ and $C$;  moreover, since any compact subgroup of $N_1$ has 
to be contained in $C$ it follows that $N_1$ is simply connected; in this case $N$ is of class $\C$.  As $d\geq 2$ it follows from 
Corollary~\ref{cor:charac} that $\Aut (N)$ is not almost algebraic. Consider on the other hand the case 
when $C$ is contained in $\overline{[N,N]}$; then $N$ not of class $\C$.   We illustrate below, through examples, that  in this case
either possibility, of being almost algebraic or otherwise, exists; the examples can be generalized to a larger class, using the theory of arithmetic groups, but here we  content ourselves noting a special case.

\subsubsection{Examples of $N$ such that $C\subset [N,N]$ and $\Aut (N)$ is  not almost algebraic} 
Let $\frak L$ be the free 2-step nilpotent Lie algebra over the $\R$-vector space $V=\R^n$, $n\geq 2$; thus $\frak L = V\oplus \wedge^2V$, and $[\frak L, \frak L]=\wedge^2V$. Let $L$ be the simply connected nilpotent Lie group associated with $\frak L$. Then $[L,L]$ is a vector group (isomorphic to $\R^{\frac 12 n(n-1)}$) and constitutes the center of $L$. Let $\Lambda$ be the lattice $\wedge^2 \Z^n$ in $ \wedge^2V=[L,L]$, and $N$ be the Lie group $L/\Lambda$. 
Then $C=[L,L]/\Lambda$ is the maximal torus in $N$ and it is also the center of $N$.  It can be seen that the group of restrictions of automorphisms of $N$ to $C$ is a countable infinite group isomorphic to $\SL(n,\Z)$.  Hence $\Aut (N)$ is not almost algebraic. 

\subsubsection{Examples of $N$ with $C$ of dimension greater that 1 for which 
$\Aut (N)$ is almost algebraic} (Note that in the light of the discussion above for this to happen $N$ has to be not of class $\C$, so $[N,N]\cap C$ is nontrivial; we shall actually have $C\subset [N,N]$). These examples, with $k\geq 2$, show in particular that the ``only if'' parts in the conclusions (1) and (2) in Corollary~\ref{cor:charac} may  not hold when the Lie group is not of class $\C$. 

{\bf a)} Let $\mathbb{H}$ be the 3-dimensional Heisenberg group, which  is a simply connected nilpotent Lie group. 
We denote by  $Z$ the  one-dimensional center of $\mathbb{H}$. Let $\Lambda$ be a lattice in $Z$. Let $k\geq 2$ and  $M=H/\Lambda $ and $N=M^k$. Let $K=Z/\Lambda$; then $C=K^k$.
Since $K$ is one-dimensional, by Theorem~\ref{characterization} $\Aut (M)$ is almost algebraic. To uphold almost algebraicity of $\Aut (N)$  we show that $\Aut (N)$ contains $(\Aut (M))^k$ as a subgroup of finite index, under the canonical identifications.  Let  $S$ be the subset of $N$ consisting of elements  whose centralizer in 
$N$ is of codimension $1$. It may be seen that $S$  consists precisely of elements of the form $(0,\dots,m_j, \dots  0)$, for some $j=1, \dots, k$,  with $m_j \notin K$.  It follows in particular that $S$ has $k$ 
connected components.   Each automorphisms of $N$ leaves the set $S$ invariant and hence permutes its connected components. A subgroup 
of finite index in $\Aut (N)$ therefore leaves invariant each of the connected components of $S$, and in turn also each copy in the cartesian 
product, viz.\ it is the subgroup  $(\Aut (M))^k$, under canonical identifications, 
as stated above. Hence $\Aut (M)$ is almost algebraic.

We also recall here, as noted in \cite{D} (p.\,451) and discussed 
in detail in \cite{PW}, that the automorpism groups $\Aut ((\SL (2,\R)\ltimes N)^k)$ is almost algebraic for all $k\in \N$.

{\bf b)} Let $\mathbb{H}$ and $Z$ be as in (a) above. 
 Let $V=\R^n$ for some $n\geq 1$ and let $W=Z\times V$, 
viewed canonically as a vector space of dimension $n+1$, as the direct sum of $Z$ and $V$. Let $z$ be a nonzero element of $Z$. Let $\Lambda$ be a lattice in $W$ with the following properties: (i) $Z+\Lambda$ is dense in $W$ and (ii) 
$z$ is not an eigenvector for an eigenvaue $\lambda$ with $|\lambda|\neq 1$ of any $\sigma \in \GL(W)$ such that $\sigma (\Lambda)=\Lambda$; it can be verified that if $\Lambda_0$ is any lattice in $W$ then the desired conditions hold for the lattice $\Lambda =g\Lambda_0$, for almost all $g\in \GL(W)$.  Now let $N=(\mathbb{H}\times V)/\Lambda$. Then $N$ is 
a nilpotent Lie group and $\overline {[N,N]}$ is the maximal compact (central) subgroup of $N$. We note that $\mathbb{H} \times V$ is a universal 
covering group of $N$, with the quotient map as a covering homomorphism. Consider any $\tau \in \Aut (N)$, and $\tilde \tau $ the 
automorphism of $\mathbb{H} \times V$ which is the lift of $\tau$. Then $\Lambda$ is $\tilde \tau$-invariant. Also, since $C$ and $[N,N]$ are  
$\tau$-invariant we get that $W$ and $Z$ are $\tilde \tau$-invariant. This means that the restriction of $\tilde \tau$ to $W$ is an automorphism, viz.
 an element of $\GL(W)$, leaving $\Lambda$ invariant and such that $z$ is an  eigenvector of $\tilde \tau$. Condition~(ii) in the choice of $\Lambda$  then implies that the restriction of $\tilde \tau$ to $Z$ is $\pm I$, $I$ being the identity map.   Since $Z+\Lambda$ 
is dense in $W$, this further implies that the restriction of $\tau$ to $C$ is either the identity automorphism or the inversion automorphism.  Thus the subgroup consisting of restrictions of automorphisms of $N$ to $C$ is finite. As $C$ is a maximal compact subgroup of $N$ and $N/C$ is simply connected 
by Theorem~\ref{characterization} we get that $\Aut(N)$ is almost algebraic. 

\begin{remark}
The above examples also show that for a connected Lie group $G$, the Lie subgroup $[G,G]$ need not be closed, and that 
condition in Theorem~\ref{characterization} needs to be formulated in terms of $\overline{[G,G]}$, rather than $[G,G]$. 
\end{remark}

\begin{remark}
Proceeding along the lines of the construction of examples in \S\,7.1.2, but replacing condition~(ii) by (ii)$' \ z$ is an eigenector with eigenvalue $\lambda$, $|\lambda|\neq 1$, of a  $\sigma \in \GL(W)$ leaving $\Lambda$ invariant, one gets a connected nilpotent Lie group $N$ such that the group of connected components $\Aut (N)/\Aut (N)^0$ has an infinite cyclic subgroup of index 2. We note however that this can happen only for $\lambda$ from a countable set, consisting of algebraic numbers. It may also be noted that when for a $\lambda$ for which such $\sigma $ and $\Lambda$ exist, the corresponding automorphism of $C=W/\Lambda$ is uniquely defined.  We shall not go into further details of this here.   

\end{remark}

\subsubsection{For a  torus $T$ not containing $C$, $F_T(N)$ need not be  almost algebraic}
By Theorem~\ref{torus-fixing} if $T$ is a torus containing $C$, $F_T(N)$ is almost algebraic whenever $N/\overline{[N,N]}C$ is simply connected. This is not true in general for tori not containing $C$. Firstly it may be observed that if $T$ is the trivial torus then $F_T(N)=\Aut (N)$ which, as we have noted, need not be almost algebraic. Examples of nontrivial $T$, not containing $C$, can be constructed using subtori of $C$ invariant under actions of subgroups of $\SL(n,\Z)$, as the group of restrictions. We illustrate this in the special case of tori invariant under the action of a unipotent element in $\SL(n,\Z)$.  

Let $u$ be a unipotent element in $\SL(n,\Z)$ of rank 2 (viz. $(u-I)^2\neq 0$ and $(u-I)^m= 0$ for some $2<m\leq n$, where $I$ and $0$ are the identity and zero matrices respectively) and $U$ be the corresponding linear transformation of $V$. Then $\wedge^2 U$ is a unipotent linear transformation of $\wedge^2V=[L,L]$, of rank at least $2$. Let  $W$ be the set of fixed points of $\wedge^2 U$. Then $W$ is a proper nonzero subspace of $[L,L]$ invariant under $\wedge^2 U$, and the factor of $\wedge^2 U$ on $[L,L]/W$ is of infinite order. 
Moreover, as $u\in \SL(n,\Z)$ it follows that $W\cap \Lambda$ is a lattice in $W$.   Let $T=W\Lambda/\Lambda$. Then $T$ is a proper nontrivial subtorus of $C$. We note that $F_T(N)$ is not almost algebraic; otherwise the subgroup  of $\Aut (N/T)$ consisting of  the factors of automorphisms of $N$ would be almost algebraic and in turn restrictions of the factors to $C/T$ would form a finite group; however this is not the case since for $\wedge^2 U^k$, $k\in \N$, the corresponding automorphisms of $C/T$ arrived at in this way are all distinct, as  the factor of $\wedge^2U$ on $[L,L]/W$ is of infinite order. It may be noted that in this case the other condition as in Theorem~\ref{torus-fixing} is satisfied, as $N/\overline {[N,N]}C =N/[N,N]$ which is a vector space. These examples point to the need for restricting the tori $T$ under consideration to those containing $C$. 

 \subsection{Automorphism groups of some semidirect products}

In the previous subsection we considered automorphism groups of nilpotent Lie groups. In this case condition (ii) of Theorem~\ref{characterization} is automatically fulfilled and as such the relevance of the condition could not be tested through these. In this subsection we discuss examples of some simple solvable groups which illustrate the need for the condition. 

Let $\mathbb{H}$ be the 3-dimensional Heisenberg group, $Z$ the center of $\mathbb{H}$, and $\Lambda$ is a lattice in  $Z$. 
Let $N=\mathbb H/\Lambda$ and  $C=Z/\Lambda$. There is a canonical action of ${\rm SL}(2,\R)$ on $N$ as a group of automorphisms, and the action  
is  trivial on $C$.
Now consider the Lie group  $G={\rm SO}(2,\R)\ltimes N$, where ${\rm SO}(2,\R)$, the special orthogonal group; it is a maximal 
compact connected one-dimensional subgroup of $\SL(2,\R)$. We note that  $G$ is a connected solvable Lie group with 
$C$ as its center; this group  is called the {\it Walnut group} (see, for example, \cite{DGS}). It is one of the simplest examples of a non-nilpotent 
Lie group which is not of class~$\C$. Now, $[G,G]=N$, so $G/\overline{[G,G]}C=G/N$ which is isomorphic to ${\rm SO}(2,\R)$. It is not simply connected and hence by Theorem~\ref{characterization} $\Aut (G)$ is not almost algebraic. We note that in this case $C$ is one-dimensional  and hence condition~(i) of the theorem is satisfied. This highlights the relevance of condition~(ii) for the theorem. We may also note that in this case  ${\rm Hom}(G/N, C)$ is 
an infinite cyclic group, so $F_C(G)/F_C(G)^0$ is infinite (cf.\ Theorem~\ref{prop-abel}) and in particular $ F_C(G)$ is not almost 
algebraic; this illustrates the need for the condition in Theorem~\ref{torus-fixing} that $G/\overline{[G,G]}T$ be simply connected.

It may also be observed that if $H$ is a noncompact closed connected subgroup of $\SL(2,\R)$ then it is either
$\SL(2,\R)$ or a simply connected solvable subgroup (including possibly a noncompact one-parameter subgroup) of $\SL(2, \R)$. In these cases for $G=H\ltimes N$ it can be readily seen that $G/\overline{[G,G]}C$ is simply connected and, since $\Aut (C)$ is finite, by Theorem~\ref{characterization}  $\Aut (G)$ is almost algebraic. 

\subsection{Examples of Lie groups with $S\cap R $ infinite}
In Corolllary~\ref{cor:charac}\,(2) we obtained a stronger result on almost algebraicity of $\Aut (G)$ under the assumption that $S\cap R$ is finite.  As mentioned earlier, while this condition  is met in most situations that one encounters in literature, it does not always hold. This is seen from the following example, which provides also an example for which the conclusion as in Corollary~\ref{cor:charac}\,(2) does not hold when $S\cap R$ is not finite. 

Let $S$ be the universal covering group of $\SL(2,\R)$. 
The center, say $D$, of $S$ is an infinite cyclic group. Let $\psi:D\to \Z$ be an isomorphism, with the latter viewed canonically as a subgroup of $\R$, 
which for notational convenience we shall denote by $L$. We now define $G$ to be the Lie group $(S\times L)/\Lambda$,  where $\Lambda$ is the 
subgroup $\{(d, \psi (d))\mid d\in D\}$. Let $\phi:S\times L \to G$ be the canonical quotient map. It is then easy to see that $\phi (L)$ is the  
radical of $G$ and $\phi (S)$ is a semisimple Levi subgroup; we note that the restrictions of $\phi$ to the copies of $S$ and $L$ in the cartesian 
product are embeddings. The intersection of $\phi (S)$ and $\phi (L)$ is 
$\phi (D\times \psi (D))$, which is an infinite cyclic group. It can be seen 
that $G$ has a one-dimensional maximal torus, which is neither contained in the radical $\phi (R)$, nor in any semisimple Levi subgroup $\phi (S)$. For this group $G$ 
itself, $\Aut (G)$ is almost algebraic, as the center does not contain any torus of positive dimension. On the other hand  consider $G'=G\times K$, 
with  $G$ as above and $K$ the circle group; we note that $K$ is then the maximal torus of the center of $G'$.  Then  $G'/\overline{[G',G']}K$ is isomorphic to $L/\phi (D)$ which is not simply  connected and hence it follows from Theorem~\ref{characterization} that $\Aut (G')$ is not almost algebraic. This shows that the statement in Corollary~\ref{cor:charac}\,(2) may not hold when $S\cap R$ is infinite. As $\Aut (G')$ is not almost algebraic and $\Aut (K)$ is finite, it follows that $F_K(G')$ is not almost algebraic. As $K$ is a maximal torus in the radical of $G'$
this shows the need for the condition of finiteness of $S\cap R$, in the last statement in Theorem~\ref{torus-fixing}. 

\medskip 
\noindent{\bf Acknowledgements} It is a pleasure to acknowledge with thanks the support received by the authors for participation in the 
International Colloquium (IC 2024), organized jointly by the Tata Institute of Fundamental Research (TIFR), Mumbai, and 
Indian Institute of Science Education and Research (IISER), Pune; the paper arose out of a discussion the authors had at the meeting. The authors would also like to record their thanks to an anonymous referee for helpful comments on an earlier version of the paper.


\begin{thebibliography}{HD}

\bibitem{BT} A. Borel and J. Tits, Groupes r\'eductifs, Inst.\ Hautes \'Etudes Sci.\ Publ.\ Math.\ 27 (1965), 55--150.


\bibitem{CS} D. Chatterjee and Riddhi Shah, Characterisation of distal actions of automorphisms on the space of
  one-parameter subgroups of Lie groups, J.\ Aust.\ Math.\ Soc.\ 119 (2025) 152--175.


\bibitem{CW} P.B. Chen and T.S. Wu, On the component group of the automorphism group of
 a Lie group, J.\ Austral.\ Math.\ Soc.\ Ser.\ A 57 (1994), 365--385.

\bibitem{C} C. Chevalley,  Th\'eorie des groupes de Lie. Tome II. Groupes alg\'ebriques (French),  Hermann \& Cie, Paris, 1951. 

\bibitem{D} S.G. Dani,  On automorphism groups of connected Lie groups, Manuscripta Math.\ 74 (1992), 445--452.


\bibitem{D-surv} S.G.  Dani, Actions of automorphism groups of Lie groups, Handbook of group actions. Vol.\ IV, 529-562, Adv.\ Lect.\ 
 Math.\ (ALM), 41, Int.\ Press, Somerville, MA, 2018. 

\bibitem{DGS} S.G. Dani, Y.\ Guivarc'h and Riddhi Shah, On the embeddability of certain infinitely divisible probability measures on Lie groups. 
Math.\ Zeit.\ 272 (2012), 361--379.

\bibitem{DM} S.G. Dani and M. McCrudden,  Convolution roots and embeddings of probability
 measures on Lie groups,  Advances in Mathematics 209 (2007), 198--211.

\bibitem{H} G. Hochschild, The structure of Lie groups. Holden-Day, Inc., San Francisco-London-Amsterdam, 1965. 

\bibitem{PW} W.H.\ Previts, and T.S. Wu, On the group of automorphisms of an analytic group, Bull.\ Austral.\,Math.\ Soc.\,64 (2001), 423--433.

\bibitem{V} V.S. Varadarajan, Lie groups, Lie Algebras, and their Representations,
Prentice-Hall, 1974.

\bibitem{W} D. Wigner, On the automorphism group of a Lie group. Proc.\ Amer.\ Math.\ Soc.\ 
 45 (1974), 140--143; Erratum: Proc.\ Amer.\ Math.\ Soc.\ 60 (1976), 376.
\end{thebibliography}
\end{document}